\documentclass{article}

\usepackage{amsmath}
\usepackage{amssymb}
\usepackage{amsthm}
\usepackage{enumerate}
\usepackage{xcolor}
\usepackage{hyperref}
\hypersetup{
    colorlinks=true,
    linkcolor=magenta,
    citecolor=cyan,
    urlcolor=magenta
}
\usepackage{amsfonts}

\newtheorem{theorem}{Theorem}[section]
\newtheorem{proposition}[theorem]{Proposition}
\newtheorem{corollary}[theorem]{Corollary}
\newtheorem{lemma}[theorem]{Lemma}
\newtheorem{rmk}[theorem]{Remark}

\newtheorem{definition}[theorem]{Definition}

\newcommand{\pp}{\mathfrak{p}}
\newcommand{\Pp}{\mathfrak{P}}

\newcommand{\Q}{\mathbb {Q}}
\newcommand{\Z}{\mathbb {Z}}

\renewcommand{\O}{\mathcal {O}}
\newcommand{\F}{\mathbb{F}}
\newcommand{\ldv}{\mathcal{L}_{\rm div}}
\newcommand{\NS}{\mathcal{N}_S}

\newcommand{\assoc}{\operatorname{assoc}}
\newcommand{\red}{\operatorname{red}}

\newcommand{\prd}{\operatorname{Prod}}

\newcommand{\qq}{\mathfrak{q}}

\title{Diophantine definability from addition and divisibility over rings of $S$-integers}
\author{ Carlos  Mart\'inez-Ranero\footnote{The author was partially supported by Proyecto VRID-Investigaci\'on No.~220.015.024-INV.}}

\begin{document}
\maketitle
\begin{abstract}
Let $K$ be a number field, let $S$ be a finite set of places of $K$ containing the archimedean places, and let $\O_{K,S}$ denote the ring of $S$-integers of $K$. We prove that multiplication is positive-existentially definable in $(\O_{K,S},0,1,+,\mid)$ if and only if the unit group $\O_{K,S}^\times$ is infinite. This extends a 1977 theorem of Lipshitz for rings of integers to rings of $S$-integers. One of the key ingredients is a uniform bound for unit equations due to Evertse, Schlickewei and Schmidt, which we use to recover multiplication on the unit group.
\end{abstract}


\section{Introduction}\label{sec:intro}

In the 1970s, Lipshitz \cite{L78} and, independently, Bel'tyukov \cite{B76} proved that one can decide algorithmically whether a system of linear equations and divisibility conditions has a solution in $\Z$. More precisely, there is an algorithm which decides, given a system
\begin{equation} \label{eqn:intro}
f_1(\mathbf{x})=0,\ \dots,\ f_a(\mathbf{x})=0,\qquad g_1(\mathbf{x})\mid h_1(\mathbf{x}),\ \dots,\ g_b(\mathbf{x})\mid h_b(\mathbf{x})
\end{equation}
in the unknowns $\mathbf{x}=(x_1,\dots,x_n)$, where the $f_i$, $g_j$ and $h_j$ are linear polynomials with integer coefficients, whether it has a solution in $\Z$. This is in sharp contrast with Hilbert's tenth problem: by the theorem of Matiyasevich \cite{Mat70}, building on work of Davis, Putnam and Robinson \cite{DPR61}, no such algorithm exists if polynomials of arbitrary degree are allowed. In the language of logic, the result of Lipshitz and Bel'tyukov says that the positive-existential theory of the structure $(\Z,0,1,+,\mid)$ is decidable.

In a series of papers, Lipshitz also studied the problem over the ring of integers $\O_K$ of a number field $K$. He showed that it is decidable when $K$ is an imaginary quadratic field \cite{L78} and undecidable when $K$ is real quadratic \cite{L78b}, and that for every number field $K$ other than $\Q$ or an imaginary quadratic field, multiplication is positive-existentially definable in $(\O_K,0,1,+,\mid)$ \cite{L77}. The latter result shows that the problem over $\O_K$ is as hard as Hilbert's tenth problem over $\O_K$, which was only recently shown to have a negative answer, by Koymans and Pagano \cite{KP24} and, independently, by Alp\"{o}ge, Bhargava, Ho and Shnidman \cite{ABHS25}. Thus, for rings of integers, the problem is decidable precisely when the unit group $\O_K^\times$ is finite, and it is natural to ask what happens for rings of $S$-integers.

So far, this question has been answered for $\Q$ and imaginary quadratic fields: Cerda-Romero and the author \cite{CM17} showed that the problem is undecidable for $\Z[1/m]$ for every integer $m\geq2$, and Hormaz\'abal and the author \cite{HM} extended this to rings of $S$-integers of imaginary quadratic fields. The problem has also been studied in other directions; see, for instance, Pheidas \cite{Phe04} for polynomial rings and Mart\'inez-Ranero, Utreras and Vidaux \cite{muv} for holomorphy subrings of global fields.

In this paper we treat arbitrary rings of $S$-integers of number fields. Let $K$ be a number field, and let $S$ be a finite set of places of $K$ containing the set $M_K^\infty$ of archimedean places. The ring of $S$-integers of $K$ is
\[
\O_{K,S}:=\{\, x\in K\colon |x|_v\leq 1 \mbox{ for all } v\notin S\,\}.
\]
Our main result is the following.

\begin{theorem}\label{mainund}
Multiplication is positive-existentially definable in the structure $(\O_{K,S},0,1,+,\mid)$ if and only if the unit group $\O_{K,S}^\times$ is infinite.
\end{theorem}

By the $S$-unit theorem, $\O_{K,S}^\times$ is infinite unless $S=M_K^\infty$ and $K$ is $\Q$ or imaginary quadratic. Thus multiplication is positive-existentially definable whenever $S\neq M_K^\infty$, and when $S=M_K^\infty$ Theorem~\ref{mainund} recovers the theorem of Lipshitz \cite{L77} by a different method. The converse direction is immediate: if $\O_{K,S}^\times$ is finite, then the positive-existential theory of $(\O_{K,S},0,1,+,\mid)$ is decidable \cite{L78,B76}, whereas Hilbert's tenth problem over $\O_{K,S}=\O_K$ has a negative answer \cite{Mat70,KP24,ABHS25}, even for equations with integer coefficients (see Theorem~\ref{thm:H10-Z}), so multiplication cannot be positive-existentially definable. The results of \cite{KP24,ABHS25}, combined with work of Eisentr\"ager \cite[Theorem~7.9]{Eis03}, give a negative answer to Hilbert's tenth problem over $\O_{K,S}$ as well, even for equations with integer coefficients (see Section~\ref{sec:proofs}). Together with Theorem~\ref{mainund}, this gives the following.

\begin{corollary}\label{undec-intro}
Suppose that the unit group $\O_{K,S}^\times$ is infinite. Then there is no algorithm which decides, given a system of the form~\eqref{eqn:intro}, whether or not it has a solution in $\O_{K,S}$. That is, the positive-existential theory of $(\O_{K,S},0,1,+,\mid)$ is undecidable.
\end{corollary}

We now describe the proof of the nontrivial direction of Theorem~\ref{mainund}, namely that multiplication is positive-existentially definable when $\O_{K,S}^\times$ is infinite. It follows the strategy of Lipshitz \cite{L77}: one first shows that the relation $x\neq0$ is positive-existentially definable, then defines multiplication on the unit group, and finally defines squaring by comparing the sizes of elements related by divisibility. However, Lipshitz measures the size of an element $a\in\O_K$ by its house, that is, the largest absolute value of its conjugates, $\max_v|\sigma_v(a)|$, and he controls it using the fact that a nonzero algebraic integer has norm of absolute value at least $1$, which is no longer true for $S$-integers. The previous extensions to rings of $S$-integers \cite{CM17,HM} deal with this failure by explicit estimates, which remain manageable because there is only one archimedean place; moreover, they define multiplication on the unit group by an argument which requires $\O_K^\times$ to be finite, and this is why they are restricted to $\Q$ and imaginary quadratic fields. The key new development is that we deal with both difficulties in a way which does not rely on these special features: we measure size by the Weil height together with the norm away from $S$, which turns divisibility into size estimates in a uniform way, and we define multiplication on the unit group using the uniform bound of Evertse, Schlickewei and Schmidt for unit equations. These two ingredients may be of independent interest.

The paper follows the three steps of this strategy. In Section~\ref{sec:nonzero} we show that the relation $x\neq0$ is positive-existentially definable, and in Section~\ref{sec:units} we define multiplication on the unit group. In Section~\ref{sec:bounds} we develop the size estimates, and in Section~\ref{sec:squaring} we use them to define squaring; multiplication then follows from the identity $(x+y)^2=x^2+2xy+y^2$. Finally, in Section~\ref{sec:proofs} we complete the proof of Theorem~\ref{mainund} and deduce Corollary~\ref{undec-intro}.

\emph{Notation.} Throughout the paper, $K$ is a number field of degree $d$ with $r_1$ real embeddings and $r_2$ pairs of complex embeddings, and $S$ is a finite set of places of $K$ containing $M_K^\infty$. We write $M_K^0$ for the set of non-archimedean places of $K$. Each $v\in M_K^0$ corresponds to a nonzero prime ideal $\pp_v$ of $\O_K$, and every nonzero prime ideal arises in this way; when $\pp=\pp_v$ we also write $|\cdot|_\pp$ for $|\cdot|_v$. We write $v_\pp$ for the normalized discrete valuation attached to a nonzero prime ideal $\pp$, and
\[
\mathcal{P}_S:=\{\pp_v\colon v\in S\cap M_K^0\}
\]
for the set of prime ideals corresponding to the non-archimedean places in $S$. By the $S$-unit theorem, $\O_{K,S}^\times$ is finitely generated of rank $r:=\#S-1$; thus $\O_{K,S}^\times$ is infinite if and only if $r\geq1$. We write $\ldv=\{=,0,1,+,\mid\}$ for the language of the structure in Theorem~\ref{mainund}.

All formulas are interpreted in $\O_{K,S}$. In particular, $a\mid b$ means that $b=ac$ for some $c\in\O_{K,S}$, and $a$ and $b$ are \emph{associated} if $a\mid b$ and $b\mid a$. A formula is \emph{positive-existential} if it is built from atomic formulas using only $\wedge$, $\vee$ and $\exists$. Subtraction and multiplication by fixed rational integers will be used freely in $\ldv$-formulas, since they can be eliminated by means of auxiliary variables and additive equations. A subset of $\O_{K,S}^n$ is \emph{positive-existentially definable} if it is defined in $\O_{K,S}$ by a positive-existential $\ldv$-formula without parameters, and, as usual, a function $\O_{K,S}^n\to\O_{K,S}$ is called positive-existentially definable if its graph is; thus multiplication is positive-existentially definable if the relation $z=xy$ is. Further notation will be introduced as it is needed.

\section{Positive-existential definability of $\O_{K,S}\setminus\{0\}$}\label{sec:nonzero}

The goal of this section is to show that the set of nonzero elements of
$\O_{K,S}$ is positive-existentially definable (see
Proposition~\ref{dif}). The idea, going back to Lipshitz~\cite{L78b,L77}, is to
consider a formula of the form
\[
\psi_{\ne}(x):=\exists \lambda,\mu,t\,\bigl(\lambda x+\mu(pt+b)=1\bigr),
\]
where $p$ is an odd prime and $b$ an integer, both to be chosen below. For
this formula to define $x\neq0$, the arithmetic progression $pt+b$ must
satisfy two requirements. On the one hand, it should contain no units, so
that $\psi_{\ne}(0)$ fails; on the other hand, for every $x\neq0$ there
should exist a term $pt+b$ of the progression such that $x$ and $pt+b$
generate the unit ideal.

Each of these requirements is easily satisfied on its own; the point is to satisfy both at once. By the Chinese remainder theorem, the second requirement holds as soon as $b$ is not divisible by $p$ (see the proof below), so the problem is to find such a progression which contains no units. For rings of integers this can be done using norms: units have norm $\pm1$, while every element of the progression has norm congruent to $b^d$ modulo $p$, so it suffices to choose $b$ with $b^d\not\equiv\pm1\pmod p$ (cf.~\cite{L77}). In $\O_{K,S}$ this argument breaks down, since units need not have norm $\pm1$. We use instead the Chebotarev density theorem to find a prime $p$ such that every unit of $\O_{K,S}$ is a square modulo some prime of $K$ above $p$; choosing $b$ to be a quadratic nonresidue modulo $p$ then ensures that the progression contains no units.

We fix some notation for the proof. For a nonzero prime ideal $\pp$ of $\O_K$, let $\O_{v_\pp}$ denote the valuation ring of $v_\pp$, $\F_\pp$ its residue field, and $\red_\pp\colon\O_{v_\pp}\to\F_\pp$ the reduction map. If $\pp\notin\mathcal{P}_S$, then $\O_{K,S}\subseteq\O_{v_\pp}$, so $\red_\pp$ is defined on $\O_{K,S}$.

\begin{proposition}\label{dif}
The set $\O_{K,S}\setminus\{0\}$ is positive-existentially
$\ldv$-definable in the structure $\O_{K,S}$.
\end{proposition}

\begin{proof}
We first choose $p$ and $b$ so that the first of the two requirements above
is satisfied.

\medskip\noindent
\textbf{Claim.} There exist an odd rational prime $p$, not below any prime
in $\mathcal{P}_S$, and an integer $0<b<p$ such that $pt+b\notin\O_{K,S}^\times$ for
every $t\in\O_{K,S}$.

\smallskip\noindent
\emph{Proof of the claim.} Let $\varepsilon_1,\dots,\varepsilon_g$ generate the finitely
generated group $\O_{K,S}^\times$ (see \cite[Ch.~I, Corollary~11.7]{Neu99}), and
let $L$ be a finite Galois extension of $\Q$ containing $K$ and
$\sqrt{\varepsilon_1},\dots,\sqrt{\varepsilon_g}$. Applying the Chebotarev density theorem
\cite[Ch.~VII, Theorem~13.4]{Neu99} to the trivial conjugacy class, we find an odd
rational prime $p$ that splits completely in $L$ and is not below any prime
in $\mathcal{P}_S$. Fix a prime $\pp$ of $K$ above $p$ and a prime $\Pp$ of $L$ above
$\pp$.

Since $p$ splits completely in $L$, we have
$\F_\Pp=\F_\pp=\F_{p}$. Moreover, since $\pp\notin\mathcal{P}_S$, each $\varepsilon_i$ is a unit
in $\O_{v_\pp}$. Since
\[
v_\Pp(\sqrt{\varepsilon_i})=\tfrac12 v_\Pp(\varepsilon_i)=0,
\]
each $\sqrt{\varepsilon_i}$ is a unit in $\O_{v_\Pp}$. Hence
$\red_\pp(\varepsilon_i)=\red_\Pp(\sqrt{\varepsilon_i})^2$ is a nonzero square in $\F_{p}$ for
every $i$, and therefore $\red_\pp(\varepsilon)$ is a nonzero square for every
$\varepsilon\in\O_{K,S}^\times$.

Now choose $0<b<p$ which is not a square modulo $p$. Since $p\in\pp$, for
every $t\in\O_{K,S}$ we have $\red_\pp(pt+b)=b$, which is not a square in
$\F_{p}$. Hence $pt+b\notin\O_{K,S}^\times$, proving the claim. $\hfill\square$

\medskip
Fix $p$ and $b$ as in the claim; it remains to check the second
requirement. In $\ldv$, the formula $\psi_{\ne}(x)$ is expressed by
\[
\exists A,B,t\,\bigl(x\mid A\wedge pt+b\mid B\wedge A+B=1\bigr),
\]
which translates to: there exist $\lambda,\mu,t\in\O_{K,S}$ such that
\[
\lambda x+\mu(pt+b)=1.
\]
We show that it defines $\O_{K,S}\setminus\{0\}$.

If $x=0$, then $\psi_{\ne}(x)$ fails: otherwise $\mu(pt+b)=1$ for some
$\mu,t\in\O_{K,S}$, so $pt+b$ would be a unit, contradicting the claim.

Suppose now that $x\neq0$. Then only finitely many maximal ideals
$\qq_1,\dots,\qq_n$ of $\O_{K,S}$ contain $x$ (if $x$ is a unit, $n=0$
and any $t$ works). For each $i$, choose a residue class $t_i$ modulo
$\qq_i$ such that $pt_i+b\notin\qq_i$. This is possible: if
$p\in\qq_i$, then $\qq_i\cap\Z=p\Z$ and $0<b<p$, so every residue class
works; otherwise, $p$ is invertible modulo $\qq_i$, and only one residue
class must be avoided. Using the Chinese remainder theorem we can find
$t\in\O_{K,S}$ such that
\[
t\equiv t_i\pmod{\qq_i}
\]
for every $i$. Then $pt+b$ lies in no maximal ideal containing $x$, so
\[
x\O_{K,S}+(pt+b)\O_{K,S}=\O_{K,S}.
\]
Hence there exist $\lambda,\mu\in\O_{K,S}$ such that
\[
\lambda x+\mu(pt+b)=1.
\]
Thus $\psi_{\ne}(x)$ holds.
\end{proof}

%%%%

\section{Multiplication on the unit group}\label{sec:units}

We next give a positive-existential definition of multiplication on the
unit group. The basic idea, which already underlies Lipshitz's
argument~\cite{L77}, is the following algebraic observation.
If $x,y,z$ are units and $z=xy$, then, for every $t$,
\[
z+ty=y(x+t),
\]
so $x+t$ and $z+ty$ are associated.

The difficulty lies in the converse: we need to know that if $x+t$ and
$z+ty$ are associated for sufficiently many values of $t$, then $z=xy$.
Lipshitz~\cite{L77} establishes a statement of this kind over $\O_K$ for
units which are sufficiently large, using properties of the norm which fail in
$\O_{K,S}$, and it seems that his techniques are not easily adapted to
our setting. We therefore take a different route: we impose the
associations above for sufficiently many translations and show that, unless
$z=xy$, each of them gives rise to a nondegenerate solution of a unit equation. A uniform
bound on the number of such solutions then gives a fixed number of
translations, independent of $x,y,z$, that forces $z=xy$.

The required uniform bound is provided by a theorem of Evertse,
Schlickewei and Schmidt~\cite[Theorem~1.1]{ESS}, which we recall next.

\begin{theorem}\label{thmESS}
Let $F$ be a field of characteristic $0$, let $n\geq 1$, let $a_1,\dots,a_n\in F^\times$ and
let $\Gamma$ be a subgroup of $(F^\times)^n$ of finite rank $\rho$. Then the number of
solutions $(u_1,\dots,u_n)\in\Gamma$ of the equation
\[
a_1u_1+\cdots+a_nu_n=1
\]
such that $\sum_{i\in I}a_iu_i\neq 0$ for every nonempty subset $I\subseteq\{1,\dots,n\}$ is
at most $\exp\bigl((6n)^{3n}(\rho+1)\bigr)$.
\end{theorem}

Solutions as in Theorem~\ref{thmESS} are called \emph{nondegenerate}. The crucial feature
for our application is that the bound is uniform in the coefficients: it depends only on the
rank of the group, so the number $N$ of associations in the formula below can be chosen
independently of $x$, $y$ and $z$.
\begin{proposition}\label{produnits}
There is a positive-existential $\ldv$-formula $\prd_u(x,y,z)$ such that,
for all $x,y,z\in\O_{K,S}$, the formula $\prd_u(x,y,z)$ holds in
$\O_{K,S}$ if and only if $x,y,z$ are units and $z=xy$.
\end{proposition}

\begin{proof}
Fix a rational prime $\ell$ not below any prime in $\mathcal{P}_S$. Then $\ell^i$ is
not a unit of $\O_{K,S}$ for any $i\geq1$. Recall that $r$ is the rank of
$\O^\times_{K,S}$. Let
\[
N:=3^{18^9(r+2)},
\]
and let $\assoc(a,b)$ denote the formula
\[
a\mid b\wedge b\mid a.
\]
Define
\[
\prd_u(x,y,z):\quad
x\mid1\ \wedge\ y\mid1\ \wedge\ z\mid1\ \wedge\
\bigwedge_{i=1}^{N}\assoc(x+\ell^i,\,z+\ell^iy).
\]
We show that $\prd_u(x,y,z)$ holds if and only if $x,y,z$ are units and $z=xy$.
Suppose first that $x,y,z$ are units and $z=xy$. Then, for every
$1\leq i\leq N$,
\[
z+\ell^iy=y(x+\ell^i),
\]
and $y$ is a unit, so $x+\ell^i$ and $z+\ell^iy$ are associated. Hence
$\prd_u(x,y,z)$ holds.

Conversely, suppose that $\prd_u(x,y,z)$ holds. Then $x,y,z$ are units.
Fix $1\leq i\leq N$. Since $x+\ell^i$ and $z+\ell^iy$ are associated, there is
a unit $\varepsilon_i$ such that
\[
x+\ell^i-\varepsilon_iz-\varepsilon_i\ell^iy=0.
\]
Dividing by $x$ gives
\begin{equation}\label{equprin3-2}
\frac{z}{x}\,\varepsilon_i-\frac{1}{x}\,\ell^i
+\frac{y}{x}\,\varepsilon_i\ell^i=1.
\end{equation}

Let
\[
\Gamma=
\{(\varepsilon,t,\varepsilon t)\colon
\varepsilon\in\O^\times_{K,S},\ t\in\langle\ell\rangle\}
\subseteq(K^\times)^3. \]
Notice that $\Gamma$
has rank $r+1$, since the projection onto the first two coordinates gives an isomorphism onto $\O^\times_{K,S}\times\langle\ell\rangle$.

By~\eqref{equprin3-2}, for each $1\leq i\leq N$ the triple
\[
(\varepsilon_i,\ell^i,\varepsilon_i\ell^i)\in\Gamma
\]
is a solution of the equation
\[
a_1u_1+a_2u_2+a_3u_3=1,\qquad\mbox{where}\quad
a_1=\frac{z}{x},\quad
a_2=-\frac{1}{x},\quad
a_3=\frac{y}{x}.
\]
Applying Theorem~\ref{thmESS} with $n=3$ and $\rho=r+1$, we see that this equation
has at most
\[
\exp(18^9(r+2))<N
\]
nondegenerate solutions in $\Gamma$.

The $N$ solutions
\[
(\varepsilon_i,\ell^i,\varepsilon_i\ell^i),\qquad 1\leq i\leq N,
\]
are pairwise distinct, since their second coordinates are distinct.
Hence at least one of them is degenerate. Fix such an $i$. Since all three terms in~\eqref{equprin3-2} are
nonzero and their sum is $1$, any vanishing proper subsum consists of
exactly two terms. There are three possibilities:
\begin{enumerate}[(i)]
\item If
\[
\frac{z}{x}\varepsilon_i-\frac{1}{x}\ell^i=0,
\]
then $\ell^i=\varepsilon_iz$, which is impossible since $\varepsilon_iz$ is a unit whereas
$\ell^i$ is not.

\item If
\[
\frac{z}{x}\varepsilon_i+\frac{y}{x}\varepsilon_i\ell^i=0,
\]
then $\ell^i=-z/y$, which is again impossible since $-z/y$ is a unit.

\item If
\[
-\frac{1}{x}\ell^i+\frac{y}{x}\varepsilon_i\ell^i=0,
\]
then $\varepsilon_i=y^{-1}$.
\end{enumerate}
Thus case (iii) must occur. Substituting $\varepsilon_i=y^{-1}$ into
\eqref{equprin3-2} gives
\[
\frac{z}{xy}=1,
\]
and hence $z=xy$. This concludes the proof of the proposition.
\end{proof}

In particular, for every fixed integer $n\geq1$, iterating $\prd_u$ shows that the relation ``$\varepsilon$ is a unit and $\eta=\varepsilon^n$'' is positive-existentially definable; we shall use such conditions freely in what follows.

%%%%
\section{Size estimates from divisibility}\label{sec:bounds}

In this section we develop the size estimates which are needed to define squaring in Section~\ref{sec:squaring}. In $\Z$, size is measured by the absolute value, and divisibility controls it: if $a\mid b$ and $b\neq0$, then $|a|\leq|b|$. In particular, a nonzero integer cannot be divisible by an integer of larger absolute value, and this is the kind of contradiction on which Section~\ref{sec:squaring} rests. Over $\O_{K,S}$ we first need a global measure of size, which takes all places of $K$ into account: Lipshitz \cite{L77} uses the house over $\O_K$, and we use the Weil height $H_K$.

The difficulty is to control $H_K$ by means of divisibility. Over $\O_K$, Lipshitz does so through the norm, which is monotone with respect to divisibility; in $\O_{K,S}$ the norm no longer has this property. We use instead the norm $\NS$ away from $S$, which is monotone with respect to divisibility but gives only partial information, and we recover the height from it by means of the product formula and suitable translations (see Lemma~\ref{lemarch}). The main result of the section is Proposition~\ref{prop:size}.

We begin by fixing the normalization of the absolute values of $K$, and then define $H_K$ and $\NS$ precisely. Let $M_K$ denote the set of places of $K$, so that $M_K=M_K^\infty\cup M_K^0$. For a nonzero prime ideal $\pp$ of $\O_K$ we write $\mathrm{N}\pp:=\#(\O_K/\pp)$ for its absolute norm. For each $v\in M_K^\infty$ we fix an embedding $\sigma_v\colon K\to\mathbb{C}$ inducing $v$, and we set $d_v=1$ or $2$ according to whether $v$ is real or complex; thus $\#M_K^\infty=r_1+r_2$ and $\sum_{v\in M_K^\infty}d_v=d$. We normalize the absolute values by
\[
|a|_\pp=\mathrm{N}\pp^{-v_\pp(a)}\quad(\pp\ne0\mbox{ prime}),\qquad |a|_v=|\sigma_v(a)|^{d_v}\quad (v\in M_K^\infty),
\]
so that the product formula $\prod_{v\in M_K}|a|_v=1$ holds for every $a\in K^\times$. Note that if $v$ is complex then $|\cdot|_v$ does not satisfy the triangle inequality; in all cases, however, we have
\begin{equation}\label{eq:weak-triangle}
|a+b|_v\leq 2^{d_v}\max\{|a|_v,|b|_v\}\qquad(v\in M_K^\infty),
\end{equation}
while $|a+b|_\pp\leq\max\{|a|_\pp,|b|_\pp\}$, with equality if $|a|_\pp\ne|b|_\pp$. Finally, we write $p_1,\dots,p_k$ for the distinct rational primes lying below the primes in $\mathcal{P}_S$, and set $P:=p_1\cdots p_k$ (so that $k=0$ and $P=1$ if $S=M_K^\infty$). If a rational integer $n$ is coprime to $P$, then $v_\pp(n)=0$ for every $\pp\in\mathcal{P}_S$.

\begin{definition}\label{def:height}
For $a\in K$, the \emph{height} of $a$ is
\[
H_K(a):=\prod_{v\in M_K}\max\{1,|a|_v\}.
\]
For $a\in\O_{K,S}\setminus\{0\}$, the \emph{norm of $a$ away from $S$} is
\[
\NS(a):=\prod_{\pp\notin\mathcal{P}_S}\mathrm{N}\pp^{v_\pp(a)}.
\]
\end{definition}

Thus $H_K$ is the (relative, multiplicative) Weil height on $K$; see \cite[Part B]{HS}. If $a\in\O_{K,S}$, then $|a|_v\leq1$ for $v\notin S$, so only the places of $S$ contribute to $H_K(a)$:
\begin{equation}\label{eq:HK-S}
H_K(a)=\prod_{v\in S}\max\{1,|a|_v\}\qquad(a\in\O_{K,S}).
\end{equation}
 The quantity $\NS(a)$ is the index of the ideal $a\O_{K,S}$ in $\O_{K,S}$, although we shall not need this interpretation. Here are the basic properties of $H_K$ and $\NS$ which we shall use.

\begin{lemma}\label{lemheights}
Let $a,b\in\O_{K,S}\setminus\{0\}$.
\begin{enumerate}
\item[(i)] We have $\NS(a)=\prod_{v\in S}|a|_v$. In particular, $1\leq\NS(a)\leq H_K(a)$.
\item[(ii)] If $a\mid b$, then $\NS(a)\leq\NS(b)$.
\item[(iii)] $H_K(a+b)\leq 2^dH_K(a)H_K(b)$, $H_K(ab)\leq H_K(a)H_K(b)$, $H_K(-a)=H_K(a)$, and $H_K(a^n)=H_K(a)^n$ for $n\geq 1$.
\item[(iv)] If $n$ is a nonzero rational integer coprime to $P$, then $\NS(n)=|n|^d$.
\end{enumerate}
\end{lemma}

\begin{proof}
(i) By the product formula,
\[
\prod_{v\in S}|a|_v=\prod_{\pp\notin\mathcal{P}_S}|a|_\pp^{-1}=\prod_{\pp\notin\mathcal{P}_S}\mathrm{N}\pp^{v_\pp(a)}=\NS(a).
\]
Since $v_\pp(a)\geq 0$ for $\pp\notin\mathcal{P}_S$, we have $\NS(a)\geq 1$, and $\NS(a)\leq H_K(a)$ follows by comparing with~\eqref{eq:HK-S} term by term.

(ii) If $b=ac$, then $c\in\O_{K,S}\setminus\{0\}$, so $\NS(b)=\NS(a)\NS(c)\geq\NS(a)$ by (i).

(iii) By~\eqref{eq:weak-triangle} and the ultrametric inequality,
\[
\max\{1,|a+b|_v\}\leq 2^{d_v}\max\{1,|a|_v\}\max\{1,|b|_v\}
\]
for every $v\in M_K$, where, in this proof only, we put $d_v=0$ for $v$ non-archimedean. Taking the product over all places gives the first inequality. The remaining statements follow from $\max\{1,st\}\leq\max\{1,s\}\max\{1,t\}$ and $\max\{1,s^n\}=\max\{1,s\}^n$ for real $s,t\geq 0$.

(iv) Since $n$ is coprime to $P$, we have $|n|_\pp=1$ for $\pp\in\mathcal{P}_S$, while $|n|_v=|n|^{d_v}$ for $v\in M_K^\infty$. The claim now follows from (i).
\end{proof}

\begin{rmk}
Writing $|a|_v=\max\{1,|a|_v\}\min\{1,|a|_v\}$ in (i) and using~\eqref{eq:HK-S}, we obtain
\begin{equation}\label{eq:NS-H}
\NS(a)=H_K(a)\prod_{v\in S}\min\{1,|a|_v\}\qquad(a\in\O_{K,S}\setminus\{0\}).
\end{equation}
Thus $\NS(a)$ is comparable to $H_K(a)$ precisely when $|a|_v$ is not small for any $v\in S$.
\end{rmk}

The first consequence of Lemma~\ref{lemheights} is a lower bound for units satisfying a congruence. In Section~\ref{sec:squaring} this is what forces $\varepsilon$ to be large.

\begin{corollary}\label{cor:large-unit}
Let $n$ be a nonzero rational integer coprime to $P$, and let $\varepsilon\in\O_{K,S}^\times$ with $\varepsilon\ne1$ and $n\mid\varepsilon-1$. Then $H_K(\varepsilon)\geq(|n|/2)^d$.
\end{corollary}

\begin{proof}
By Lemma~\ref{lemheights}, $|n|^d=\NS(n)\leq\NS(\varepsilon-1)\leq H_K(\varepsilon-1)\leq 2^dH_K(\varepsilon)$.
\end{proof}

We now prepare the proof of Proposition~\ref{prop:size}. By~\eqref{eq:NS-H}, the task is to find a translate $a+j$, or $a+j\eta$ with $\eta$ a unit, which is not small at any place of $S$. At a prime $\pp\in\mathcal{P}_S$ this is controlled by valuations: if $v_\pp(a)\ne v_\pp(b)$, then $|a+b|_\pp=\max\{|a|_\pp,|b|_\pp\}$. At an archimedean place $v$, the number $\sigma_v(a)+j$ may be very small for a particular integer $j$, but this can happen for at most one $j$; hence among any $r_1+r_2+1$ integers there is one for which $\sigma_v(a)+j$ is not small at any archimedean place. This motivates the following definition.

\begin{definition}\label{def:admissible}
A set $J$ of positive integers is \emph{admissible} if $\#J=r_1+r_2+1$ and every $j\in J$ is coprime to $P$. For such $J$ we set
\[
c_J:=(2\max J+1)^{-1}.
\]
\end{definition}

\begin{lemma}\label{lemarch}
Let $J$ be admissible, and for each $v\in M_K^\infty$ let $z_v$ be a complex number.
Then there exists $j\in J$ such that
\[
|z_v+j|\geq c_J\max\{1,|z_v|\}
\quad\text{for every }v\in M_K^\infty.
\]
\end{lemma}

\begin{proof}
For each $v\in M_K^\infty$, there is at most one integer $i$
such that $|z_v+i|<\frac12$. Since $J$ contains more elements
than there are archimedean places, we may choose $j\in J$
such that $|z_v+j|\geq\frac12$ for every $v\in M_K^\infty$.

Fix such a $j$, and write $j^*:=\max J$. For every $v\in M_K^\infty$ we have, on the one hand,
\[
1\leq 2|z_v+j|\leq(2j^*+1)|z_v+j|,
\]
and on the other hand the triangle inequality gives
\[
|z_v|\leq |z_v+j|+j
\leq |z_v+j|+2j^*|z_v+j|
=(2j^*+1)|z_v+j|.
\]
Combining these inequalities yields the desired conclusion.
\end{proof}

Here is the main result of the section.

\begin{proposition}\label{prop:size}
Let $J$ be admissible, and let $a\in\O_{K,S}\setminus\{0\}$.
\begin{enumerate}
\item[(i)] Suppose that $v_\pp(a)\ne0$ for every $\pp\in\mathcal{P}_S$. If $b\in\O_{K,S}\setminus\{0\}$ and $a+j\mid b$ for every $j\in J$, then
\[
H_K(a)\leq c_J^{-d}\,\NS(b).
\]
\item[(ii)] Let $\eta\in\O_{K,S}^\times$, and suppose that $v_\pp(a)\ne v_\pp(\eta)$ for every $\pp\in\mathcal{P}_S$. Then there exists $j\in J$ such that $a+j\eta\ne0$ and
\[
\NS(a+j\eta)\geq c_J^d\,\frac{H_K(\eta)}{H_K(a)}.
\]
\end{enumerate}
\end{proposition}

\begin{proof}
(i) Applying Lemma~\ref{lemarch} with $z_v=\sigma_v(a)$, we obtain $j\in J$ such that $|\sigma_v(a)+j|\geq c_J\max\{1,|\sigma_v(a)|\}$ for every $v\in M_K^\infty$. Raising this to the power $d_v$ gives
\[
|a+j|_v\geq c_J^{d_v}\max\{1,|a|_v\}\qquad(v\in M_K^\infty);
\]
in particular, $a+j\ne0$. For $\pp\in\mathcal{P}_S$ we have $v_\pp(j)=0\ne v_\pp(a)$, and hence $|a+j|_\pp=\max\{1,|a|_\pp\}$. Taking the product over $v\in S$ and using~\eqref{eq:HK-S} and Lemma~\ref{lemheights} (i) and (ii), we obtain
\[
c_J^dH_K(a)\leq\prod_{v\in S}|a+j|_v=\NS(a+j)\leq\NS(b).
\]

(ii) Since $\eta$ is a unit, $\sigma_v(\eta)\ne0$ for every $v\in M_K^\infty$. Applying Lemma~\ref{lemarch} with $z_v=\sigma_v(a)/\sigma_v(\eta)$ and multiplying the resulting inequalities by $|\sigma_v(\eta)|$, we obtain $j\in J$ such that
\[
|a+j\eta|_v\geq c_J^{d_v}\max\{|a|_v,|\eta|_v\}\qquad(v\in M_K^\infty).
\]
In particular, $a+j\eta\ne0$. For $\pp\in\mathcal{P}_S$ we have $v_\pp(j\eta)=v_\pp(\eta)\ne v_\pp(a)$, and hence $|a+j\eta|_\pp=\max\{|a|_\pp,|\eta|_\pp\}$. Using Lemma~\ref{lemheights} (i) and the elementary inequality $\max\{s,t\}\geq\max\{1,t\}\min\{1,s\}$, valid for real $s,t\geq0$, we get
\[
\NS(a+j\eta)\geq c_J^d\prod_{v\in S}\max\{|a|_v,|\eta|_v\}\geq c_J^d\prod_{v\in S}\max\{1,|\eta|_v\}\prod_{v\in S}\min\{1,|a|_v\}.
\]
By~\eqref{eq:HK-S} and~\eqref{eq:NS-H}, the two products on the right are $H_K(\eta)$ and $\NS(a)/H_K(a)$, respectively. Since $\NS(a)\geq1$ by Lemma~\ref{lemheights} (i), we conclude that
\[
\NS(a+j\eta)\geq c_J^d\,\frac{H_K(\eta)}{H_K(a)}.
\]
\end{proof}

\begin{rmk}
The hypotheses on $v_\pp(a)$, $\pp\in\mathcal{P}_S$, are not automatic, and they are the reason for the rescaling of $x$ by $p^\alpha$ and of $y$ by $p^{2\alpha}$ in the next section.
\end{rmk}

We record for future use the following elementary analogue of the Euler--Fermat theorem.

\begin{lemma}\label{lemunitinforder}
Let $\xi\in\O_{K,S}^\times$ be a unit of infinite order, and let $a\in\O_{K,S}\setminus\{0\}$. Then there exists an integer $n\geq1$ such that $\varepsilon:=\xi^n$ is a unit of infinite order with $a\mid \varepsilon-1$.
\end{lemma}

\begin{proof}
Since $\O_{K,S}/a\O_{K,S}$ is finite, the image of $\xi$ in its unit group has finite order, say $n\geq1$. Then $\varepsilon:=\xi^n$ satisfies $a\mid\varepsilon-1$, and it has infinite order since $\xi$ does.
\end{proof}

%%%%
\section{Definability of the squaring function}\label{sec:squaring}

We now come to the main step of the proof: we construct a positive-existential formula $\varphi_{\rm sq}(x,y)$ which holds in $\O_{K,S}$ if and only if $y=x^2$, so that it defines the graph of the squaring function. We first write down the formula, then explain the role of each of its parts, and finally prove that it has this property (Theorem~\ref{mainthm}).

The formula involves some auxiliary integers and constants, depending only on $K$ and $S$, which we now fix (see also Remark~\ref{rmk:parameters}). Recall from Section~\ref{sec:bounds} that $p_1,\dots,p_k$ are the rational primes below the primes in $\mathcal{P}_S$ and that $P=p_1\cdots p_k$.
\begin{itemize}
\item Set $m:=r_1+r_2+1$ and $j_s:=1+Ps$ for $0\leq s<3m$. We split these $3m$ integers into three consecutive sets of $m$ elements: for $i=0,1,2$, let
\[
J_i:=\{j_{im},j_{im+1},\dots,j_{(i+1)m-1}\}.
\]
\item Set $j_{\max}:=1+P(3m-1)$, the largest of the $j_s$, and set $c:=(2j_{\max}+1)^{-1}$ and $C:=(2/c)^d$.
\item For a multi-index $\alpha=(\alpha_1,\dots,\alpha_k)\in\Z_{\geq0}^k$ we write $p^\alpha:=\prod_{i=1}^kp_i^{\alpha_i}$, and we set $\mathcal{A}:=\{0,1,\dots,3d\}^k$. (If $S=M_K^\infty$, then $\mathcal{A}$ consists only of the empty multi-index, and $p^\alpha=1$.)
\item Let $q$ be a rational prime, coprime to $P$, with $q>64c^{-5}$.
\end{itemize}
Each $j_s$ is congruent to $1$ modulo every $p_i$, and hence is coprime to $P$. Thus each $J_i$ is admissible in the sense of Definition~\ref{def:admissible}, the sets $J_0$, $J_1$ and $J_2$ are pairwise disjoint, and $c_{J_i}\geq c$ for $i=0,1,2$. They will play the role of the admissible sets of translations in Proposition~\ref{prop:size}.

We can now write down the formula. Let $\varphi_{\rm sq}(x,y)$ be
\begin{equation}\label{sqformula}
(x=0\wedge y=0)\vee \left(\psi_{\neq}(x)\wedge\psi_{\neq}(y)\wedge \exists \varepsilon,\eta,\delta\ \psi(x,y,\varepsilon,\eta,\delta)\right),
\end{equation}
where $\psi_{\ne}$ is the formula of Proposition~\ref{dif}, and $\psi(x,y,\varepsilon,\eta,\delta)$ is the conjunction of the following four formulas:
\begin{subequations}\label{eq:sq-conditions}
\begin{align}
\psi_1(\varepsilon):\quad &\varepsilon\mid 1\wedge \psi_{\neq}(\varepsilon-1)\wedge q\mid\varepsilon-1, \label{eq:sq-psi1}\\
\psi_2(\varepsilon,\eta,\delta):\quad &\eta=\varepsilon^5\wedge \delta=\varepsilon^{10}, \label{eq:sq-psi2}\\
\psi_3(x,y,\varepsilon):\quad &\bigwedge_{\alpha\in\mathcal{A}}\ \bigvee_{i=0}^{2}\ \bigwedge_{j\in J_i}
\left(p^{\alpha}x+j\mid \varepsilon-1\wedge p^{2\alpha}y+j\mid \varepsilon-1\right), \label{eq:sq-psi3}\\
\psi_4(x,y,\eta,\delta):\quad &\bigwedge_{\alpha\in\mathcal{A}}\ \bigwedge_{j\in J_0}
 p^{\alpha}x+j\eta\mid p^{2\alpha}y-j^2\delta. \label{eq:sq-psi4}
\end{align}
\end{subequations}

All the index sets and constants involved are fixed; in particular, $p^{\alpha}$, $p^{2\alpha}$, $j$ and $j^2$ are rational integers, and the conditions $\eta=\varepsilon^5$ and $\delta=\varepsilon^{10}$ can be expressed as explained at the end of Section~\ref{sec:units}. Thus $\varphi_{\rm sq}$ is a positive-existential $\ldv$-formula.

Let us explain the role of the various parts of $\varphi_{\rm sq}$. The idea is that, if $y\neq x^2$, the nonzero element $p^{2\alpha}(y-x^2)$ is forced to have a divisor which is too large. Condition $\psi_1$ gives a lower bound for the height of the unit $\varepsilon$: since $\varepsilon\neq1$ and $q\mid\varepsilon-1$, Corollary~\ref{cor:large-unit} shows that $\varepsilon$ is large, independently of $x$ and $y$. Condition $\psi_3$ gives upper bounds for the heights of $p^\alpha x$ and $p^{2\alpha}y$ in terms of that of $\varepsilon$, by Proposition~\ref{prop:size} (i). Condition $\psi_2$ provides the powers $\eta=\varepsilon^5$ and $\delta=\eta^2$. Finally, since $p^\alpha x+j\eta$ always divides $p^{2\alpha}x^2-j^2\delta$, condition $\psi_4$ implies that $p^\alpha x+j\eta$ divides $p^{2\alpha}(y-x^2)$, and by Proposition~\ref{prop:size} (ii) and the previous bounds, this divisor is too large for some $j$ unless $y=x^2$. With these estimates, the exponent $5$ is the smallest one for which this comparison works.

\begin{theorem}\label{mainthm}
Suppose that $r\geq1$. Then for all $x,y\in\O_{K,S}$, the formula $\varphi_{\rm sq}(x,y)$ holds in $\O_{K,S}$ if and only if $y=x^2$. In particular, the graph of the squaring function is positive-existentially definable in the $\ldv$-structure $\O_{K,S}$.
\end{theorem}

\begin{proof}
First suppose that $y=x^2$. If $x=0$, then the first disjunct in~\eqref{sqformula} holds. We may therefore assume that $x\ne0$, and hence $y\ne0$; thus $\psi_{\neq}(x)$ and $\psi_{\neq}(y)$ hold by Proposition~\ref{dif}. For each $\alpha\in\mathcal{A}$ there is at most one $j$ such that $p^\alpha x+j=0$ and at most one $j$ such that $p^{2\alpha}y+j=0$. Since the sets $J_0$, $J_1$ and $J_2$ are pairwise disjoint, there is $i(\alpha)\in\{0,1,2\}$ such that none of the elements $p^\alpha x+j$ and $p^{2\alpha}y+j$, with $j\in J_{i(\alpha)}$, is zero. Since $r\geq1$, there is a unit $\xi$ of infinite order; this is the only point in the proof where the hypothesis $r\geq1$ is used. Applying Lemma~\ref{lemunitinforder} to $\xi$ and to the product of $q$ and all these elements (for all $\alpha\in\mathcal{A}$), we obtain a unit $\varepsilon$ of infinite order such that $q\mid\varepsilon-1$ and~\eqref{eq:sq-psi3} holds. Since $\varepsilon\ne 1$, \eqref{eq:sq-psi1} holds by Proposition~\ref{dif}. Let $\eta=\varepsilon^5$ and $\delta=\varepsilon^{10}$, so that~\eqref{eq:sq-psi2} holds by Proposition~\ref{produnits}. Finally, \eqref{eq:sq-psi4} holds because
\[
p^{2\alpha}y-j^2\varepsilon^{10}=(p^\alpha x+j\varepsilon^5)(p^\alpha x-j\varepsilon^5).
\]
Thus all four conditions in~\eqref{eq:sq-conditions} are satisfied, and hence $\varphi_{\rm sq}(x,y)$ holds.

We now prove the converse. If the first disjunct of~\eqref{sqformula} holds, there is nothing to prove, so assume that the second disjunct holds. In particular, $x\ne0$ and $y\ne0$. Fix $\varepsilon$, $\eta$ and $\delta$ such that $\psi(x,y,\varepsilon,\eta,\delta)$ holds. Then $\varepsilon$ is a unit, $\varepsilon\ne1$, $q\mid\varepsilon-1$, $\eta=\varepsilon^5$ and $\delta=\varepsilon^{10}$. We first choose the rescaling.

\medskip
\noindent\textbf{Claim A.} There exists a multi-index $\alpha\in\mathcal{A}$ such that $v_\pp(p^\alpha x)\ne0$, $v_\pp(p^{2\alpha}y)\ne0$ and $v_\pp(p^\alpha x)\ne5v_\pp(\varepsilon)$ for every $\pp\in\mathcal{P}_S$.

\noindent\textit{Proof of Claim A.}
For each $1\leq i\leq k$, choose
\begin{align*}
\alpha_i\in \{0,\dots,3d\}\setminus\Bigg(&\left\{\frac{-v_\pp(x)}{e(\pp:p_i)}\colon \pp\mid p_i\right\}
\cup\left\{\frac{-v_\pp(y)}{2e(\pp:p_i)}\colon \pp\mid p_i\right\}\\
&\cup\left\{\frac{5v_\pp(\varepsilon)-v_\pp(x)}{e(\pp:p_i)}\colon \pp\mid p_i\right\}\Bigg)
\end{align*}
arbitrarily, where $\pp$ ranges over the primes in $\mathcal{P}_S$ above $p_i$ and $e(\pp:p_i)$ denotes the ramification index. This is possible since there are at most $d$ primes above $p_i$, so at most $3d$ values have to be avoided. Since $v_\pp(p^\alpha x)=\alpha_i e(\pp:p_i)+v_\pp(x)$ and $v_\pp(p^{2\alpha}y)=2\alpha_i e(\pp:p_i)+v_\pp(y)$ for $\pp\mid p_i$, the multi-index $\alpha=(\alpha_1,\dots,\alpha_k)$ is as required. $\hfill\square$

\medskip
Fix such an $\alpha$, and write
\[
X:=p^\alpha x,\qquad Y:=p^{2\alpha}y.
\]

\medskip
\noindent\textbf{Claim B.} We have
\begin{equation}\label{eq:sq-claimB}
H_K(X)\leq CH_K(\varepsilon)\quad\mbox{and}\quad H_K(Y)\leq CH_K(\varepsilon).
\end{equation}

\noindent\textit{Proof of Claim B.} Since~\eqref{eq:sq-psi3} holds, there is $i\in\{0,1,2\}$ such that $X+j\mid\varepsilon-1$ and $Y+j\mid\varepsilon-1$ for every $j\in J_i$. By Claim A, $v_\pp(X)\ne0$ and $v_\pp(Y)\ne0$ for every $\pp\in\mathcal{P}_S$. Applying Proposition~\ref{prop:size} (i) with $J=J_i$ to the nonzero element $\varepsilon-1$, and using Lemma~\ref{lemheights} and $c_{J_i}\geq c$, we obtain
\[
H_K(X)\leq c^{-d}\NS(\varepsilon-1)\leq c^{-d}H_K(\varepsilon-1)\leq(2/c)^dH_K(\varepsilon)=CH_K(\varepsilon),
\]
and similarly for $Y$. $\hfill\square$

\medskip
Suppose, towards a contradiction, that $D:=Y-X^2\ne0$. For every $j\in J_0$, condition~\eqref{eq:sq-psi4} gives
\[
X+j\varepsilon^5\mid Y-j^2\varepsilon^{10},
\]
while
\[
X+j\varepsilon^5\mid X^2-j^2\varepsilon^{10}.
\]
Hence $X+j\varepsilon^5\mid D$ for every $j\in J_0$. By Claim A, we may apply Proposition~\ref{prop:size} (ii) with $J=J_0$, $a=X$ and $\eta=\varepsilon^5$, and so, using $c_{J_0}\geq c$ and $H_K(\varepsilon^5)=H_K(\varepsilon)^5$ (Lemma~\ref{lemheights} (iii)), there is $j\in J_0$ such that $X+j\varepsilon^5\ne0$ and
\[
\NS(X+j\varepsilon^5)\geq c^d\frac{H_K(\varepsilon)^5}{H_K(X)}.
\]
Combining this with Claim B, Lemma~\ref{lemheights} and the divisibility $X+j\varepsilon^5\mid D$, we obtain
\begin{align}
 c^dC^{-1}H_K(\varepsilon)^4
 &\leq c^d\frac{H_K(\varepsilon)^5}{H_K(X)}
 \leq \NS(X+j\varepsilon^5) \notag\\
 &\leq \NS(D)
 \leq H_K(D)
 \leq 2^dH_K(Y)H_K(X)^2 \notag\\
 &\leq 2^dC^3H_K(\varepsilon)^3. \label{eq:sq-final-chain}
\end{align}
Canceling $H_K(\varepsilon)^3$ in~\eqref{eq:sq-final-chain} gives
\begin{equation}\label{eq:sq-upper-epsilon}
H_K(\varepsilon)\leq 2^dc^{-d}C^4=(32c^{-5})^d.
\end{equation}
On the other hand, by $\psi_1$, Corollary~\ref{cor:large-unit} and our choice of $q$,
\begin{equation}\label{eq:sq-lower-epsilon}
H_K(\varepsilon)\geq(q/2)^d>(32c^{-5})^d,
\end{equation}
contradicting~\eqref{eq:sq-upper-epsilon}. Therefore $Y=X^2$, that is, $p^{2\alpha}(y-x^2)=0$, and so $y=x^2$. This concludes the proof.
\end{proof}

\section{Undecidability}\label{sec:proofs}

We first note that Theorem~\ref{mainund} follows from Theorem~\ref{mainthm}. Indeed, suppose that $r\geq1$. Since
\[
2xy=(x+y)^2-x^2-y^2
\]
and $\O_{K,S}$ has characteristic $0$, multiplication is positive-existentially definable from squaring. The converse direction was explained after the statement of Theorem~\ref{mainund}.

\begin{rmk}\label{rmk:parameters}
The resulting formula for multiplication contains no parameters: the only constants occurring in it, and in the formulas $\psi_{\ne}$, $\prd_u$ and $\varphi_{\rm sq}$ on which it is built, are rational integers, which can be expressed in $\ldv$ by means of $1$ and $+$. It is, however, not uniform in $K$ and $S$: the integers $p$ and $b$ of Proposition~\ref{dif}, the prime $\ell$ and the number $N$ of Proposition~\ref{produnits}, and the prime $q$, the sets $J_i$ and the multi-indices in $\mathcal{A}$ of Section~\ref{sec:squaring} all depend on $K$ and $S$.
\end{rmk}

The rest of this section is devoted to the proof of Corollary~\ref{undec-intro}. The arguments are standard, and probably folklore; we include them for completeness, and in order to make precise that undecidability holds already for equations with integer coefficients, which is what is needed for the language $\ldv$.

The undecidability results quoted above concern polynomial equations whose coefficients may lie in $\O_{K,S}$. These constants cannot simply be named in $\ldv$: in general they are not even definable without parameters, since an automorphism of $K$ preserving $\O_{K,S}$ preserves addition and divisibility (for instance, $\sqrt2$ and $-\sqrt2$ satisfy the same parameter-free formulas in $\Z[\sqrt2]$). The following standard observation allows us to remove them.

\begin{lemma}\label{lemcoef}
Let $R$ be a subring of $K$ which is finitely generated as a $\Z$-algebra and has fraction field $K$. Then for every system of polynomial equations with coefficients in $R$ one can effectively construct a system of polynomial equations with coefficients in $\Z$ which has a solution over $R$ if and only if the original system does.
\end{lemma}
\begin{proof}
Fix a finite presentation
\[
R\simeq\Z[X_1,\dots,X_h]/(g_1,\dots,g_e),
\]
and let $\theta_1,\dots,\theta_h$ be the images of $X_1,\dots,X_h$. Replace each coefficient of the given system by a polynomial expression in $\theta_1,\dots,\theta_h$ with integer coefficients, replace $\theta_1,\dots,\theta_h$ by new variables $t_1,\dots,t_h$, and add the equations $g_\lambda(t_1,\dots,t_h)=0$ for $\lambda=1,\dots,e$.

Clearly, if the original system has a solution, then so does the new one (take $t_i=\theta_i$). Conversely, a solution $t_1,\dots,t_h\in R$ of the equations $g_\lambda=0$ determines a homomorphism of $\Z$-algebras $\phi\colon R\rightarrow R$ with $\phi(\theta_i)=t_i$. We claim that $\phi$ is injective. Indeed, a nonzero element $\kappa$ of its kernel is algebraic over $\Q$, so it satisfies an equation with integer coefficients and nonzero constant term $n_0$; it follows that $n_0\in\ker\phi$, which is impossible since $\phi(n_0)=n_0\ne 0$. Therefore, $\phi$ extends to a field embedding $K\rightarrow K$, which is an automorphism of finite order, say $n$. Since $\phi(R)\subseteq R$, we have $R=\phi^n(R)\subseteq\phi(R)\subseteq R$, so $\phi(R)=R$. Applying $\phi^{-1}$ to a solution of the new system we obtain a solution of the original one.
\end{proof}

\begin{theorem}\label{thm:H10-Z}
There is no algorithm which decides, given a polynomial $f\in\Z[x_1,\dots,x_n]$, whether $f$ has a zero in $\O_{K,S}^n$.
\end{theorem}
\begin{proof}
Koymans and Pagano \cite{KP24} and, independently, Alp\"oge, Bhargava, Ho and Shnidman \cite{ABHS25} proved that Hilbert's tenth problem has a negative answer over $\O_K$, and Eisentr\"ager \cite[Theorem~7.9]{Eis03} showed that $\O_K$ is Diophantine over $\O_{K,S}$; hence Hilbert's tenth problem over $\O_{K,S}$, for polynomials with coefficients in $\O_{K,S}$, also has a negative answer. To apply Lemma~\ref{lemcoef}, we note that $\O_{K,S}$ is finitely generated as a $\Z$-algebra: if $h_K$ denotes the class number of $K$, $a_\pp\in\O_K$ is a generator of $\pp^{h_K}$ for each $\pp\in\mathcal{P}_S$, and $a:=\prod_{\pp\in\mathcal{P}_S}a_\pp$, then $\O_{K,S}=\O_K[a^{-1}]$. Lemma~\ref{lemcoef} now allows us to replace the coefficients by integers.
\end{proof}

\begin{proof}[Proof of Corollary~\ref{undec-intro}]
Suppose that $r\geq1$. By Theorem~\ref{mainund}, every polynomial equation with integer coefficients can be effectively translated into a positive-existential $\ldv$-sentence, which in turn can be effectively written as a finite disjunction of sentences, each asserting that a system of the form~\eqref{eqn:intro} has a solution. An algorithm deciding the solvability of systems of the form~\eqref{eqn:intro} would therefore contradict Theorem~\ref{thm:H10-Z}.
\end{proof}

\subsection*{Statement on the use of AI tools}
During the preparation of this manuscript the author used Claude (Anthropic), an AI assistant, for proofreading, including checking the computations of the constants and the hypotheses of the statements, for advice on the structure of the paper, and for grammar and style, including drafting parts of the expository text. The mathematical results and the strategy of proof are due to the author, who reviewed and edited all content produced with the assistance of this tool and takes full responsibility for the content of the paper. This statement was itself drafted with the assistance of the same tool.

\noindent Carlos Mart\'inez-Ranero\\
Email: cmartinezr@udec.cl\\
Homepage: www2.udec.cl/~cmartinezr\\

\noindent Address: \\
Universidad de Concepci\'on, Concepci\'on, Chile\\
Facultad de Ciencias F\'isicas y Matem\'aticas\\
Departamento de Matem\'atica\\

\end{document}